\documentclass{ifacconf}

\usepackage{graphicx}      
\usepackage{natbib}        

\usepackage[ansinew]{inputenc}
\usepackage{amsmath, amsfonts}
\usepackage{amssymb}
\usepackage{subfigure}
\usepackage{multirow}
\newtheorem{proposition}{Proposition}
\newtheorem{corollary}{Corollary}
\theoremstyle{definition}
\newtheorem{definition}{Definition}
\newtheorem{remark}{Remark}

\usepackage{verbatim}
\usepackage{array}
\usepackage{float}
\usepackage{subfigure}
\usepackage{algorithm, algorithmicx, algpseudocode}

\newcommand{\guillemets}[1]{``#1''}
\newcommand{\ve}[1]{\mathbf{#1}}

\begin{document}
\begin{frontmatter}

\title{Koopman operator framework for spectral analysis and identification of infinite-dimensional systems} 


\author[First]{A. Mauroy}

\address[First]{Department of Mathematics and Namur Institute for Complex Systems (naXys), University of Namur, Belgium \\(e-mail: alexandre.mauroy@unamur.be).}

\begin{abstract} 
We consider Koopman operator theory in the context of nonlinear infinite-dimensional systems, where the operator is defined over a space of bounded continuous functionals. The properties of the Koopman semigroup are described and a finite-dimensional projection of the semigroup is proposed, which provides a linear finite-dimensional approximation of the underlying infinite-dimensional dynamics. This approximation is used to obtain spectral properties from data, a method which can be seen as a generalization of Extended Dynamic Mode Decomposition for infinite-dimensional systems. Finally, we exploit the proposed framework to identify (a finite-dimensional approximation of) the Lie generator associated with the Koopman semigroup. This approach yields a linear method for nonlinear PDE identification, which is complemented with theoretical convergence results.
\end{abstract}


\end{frontmatter}

\section{Introduction}

Koopman operator theory is a powerful framework which provides an alternative approach to dynamical systems. Through the so-called Koopman (or composition) operator (\cite{Koopman}), nonlinear dynamical systems are approximated by linear, higher-dimensional systems that are amenable to systematic analysis. Under the impulse of the seminal work \cite{Mezic}, the Koopman operator framework has grown in popularity over the last decade in dynamical systems theory, and more recently has attracted attention in control theory (see \cite{MSM_book} and references therein). However, the research effort has focused on finite-dimensional systems, and little work has been devoted to infinite-dimensional dynamical systems, such as nonlinear partial differential equations (PDEs). In this context, one can mention the early work by \cite{Banks}, where (an equivalent of) the Koopman operator was defined on a separable $L^2$ space of functionals (themselves defined on a Hilbert space). Later on, \cite{Neuberger} followed a different path, defining the composition operator in the space of bounded continuous functionals and investigating the properties of the associated Lie generator. In the recent work by \cite{Farkas}, this approach was pushed further and leveraged in the framework of strongly continuous semigroups by considering a space of bounded continuous functionals equipped with a mixed topology. Finally, it is also recently that the Koopman operator framework has been considered to study the spectral properties of nonlinear PDEs (\cite{Mezic_book,Nakao_PDE}).

In this paper, we further exploit and investigate the Koopman operator framework for systems described by infinite-dimensional differential equations. We do not rely on strongly continuous semigroup theory, but rather adopt the approach proposed by \cite{Neuberger}. In this context, the Lie generator is shown to be related to a Gâteaux derivative and a finite-dimensional approximation of the Koopman operator is presented, which allows to approximate a nonlinear, infinite-dimensional system by a linear, finite-dimensional one. Based on this framework, our main contributions are twofold. First, we complement the spectral analysis in \cite{Nakao_PDE} by developing a data-driven method to compute the spectral properties of the Koopman operator. This method is a generalization of the Extended Dynamic Mode Decomposition (EDMD) method (\cite{Rowley_EDMD}) to infinite-dimensional systems and is more general than the mere application of standard EDMD to spatially discretized PDEs. In fact, while the method is developed for general basis functionals, it narrows down to the EDMD method only for a specific choice of basis functionals. Second, combining our previous work (\cite{Mauroy_Goncalves2}) with the above framework, we propose a novel method for nonlinear identification of infinite-dimensional systems. This method relies on a linear estimation of the Koopman generator, similarly to \cite{Kaiser, Klus, Klus_other} in the context of finite-dimensional (possibly stochastic) systems, but in contrast our proposed method does not require to evaluate time derivatives. In this sense, it is an indirect alternative method to recent direct methods for data-driven discovery of nonlinear PDEs (\cite{Kutz_PDE, PDE_ident_neural_net, Yuan_PDE, Grigoriev}).

The rest of the paper is organized as follows. In Section \ref{sec:theory}, we introduce the Koopman operator framework for infinite-dimensional systems, with a focus on the Lie generator and finite-dimensional approximation. Section \ref{sec:spectral} is devoted to spectral analysis and presents the generalized Extended Dynamic Mode Decomposition (EDMD) method. The identification method for infinite-dimensional systems is presented in Section \ref{sec:identif} along with theoretical convergence results and numerical examples. Finally, concluding remarks are given in Section \ref{sec:conclu}. 

\section{Koopman operator theory for infinite-dimensional systems}
\label{sec:theory}

\subsection{Koopman semigroup}

We consider (infinite-dimensional) dynamical systems of the form
\begin{equation}
\label{eq:abstract_diff}
\dot{u} = W(u) \,, \qquad u \in \mathcal{U}
\end{equation}
where $\mathcal{U}$ is a separable Hilbert space and $W:\mathcal{D}(W) \to \mathcal{U}$ is a nonlinear operator, with $\mathcal{D}(W)$ the domain of $W$. If the system is described by a PDE, then $W$ is typically a differential operator. Moreover, we assume that $W$ generates a (possibly nonlinear) semiflow $(\varphi_t)_{t\geq 0}: \mathcal{U} \to \mathcal{U}$, i.e. $u(t)=\varphi^t(u)$ is a classical solution to the abstract differential equation \eqref{eq:abstract_diff} associated with the initial condition $u_0 \in \mathcal{D}(W)$. We also make the standing assumption that each map $\varphi^t:\mathcal{U} \to \mathcal{U}$, $t\geq 0$, is continuous and that the mapping $t \mapsto \varphi^t(u)$ is continuous from $\mathbb{R}^+$ into $\mathcal{U}$ (strong continuity).

The semigroup of Koopman operators (or Koopman semigroup in short) associated with \eqref{eq:abstract_diff} is defined on a space of \guillemets{observable-functionals.}
\begin{definition}[Koopman semigroup] Consider the space $\mathcal{E}$ of complex-valued functionals $\zeta: \overline{\mathcal{U}} \to \mathbb{C}$, where the definition domain $\overline{\mathcal{U}} \subset \mathcal{D}(W)$ is invariant under $\varphi^t$. The semigroup of Koopman operators $(K^t)_{t \geq 0}$ associated with the semiflow $(\varphi_t)_{t\geq 0}$ is defined by $K^t \zeta = \zeta \circ \varphi^t$, $\zeta \in \mathcal{E}$. \hfill $\diamond$
\end{definition}
In the rest of the paper, $\mathcal{E}$ is the space of bounded continuous functionals, endowed with the supremum norm $\|\zeta\|=\sup_{u \in \overline{\mathcal{U}}} |\zeta(u)|$, i.e. $\mathcal{E}=C(\overline{\mathcal{U}})$.

\subsection{Lie generator}

Following the work by \cite{Neuberger}, we define the Lie generator of the Koopman semigroup.
\begin{definition}[Lie generator]
\label{def:Lie_gen}
The Lie generator of the semigroup $(K^t)_{t \geq 0}$ is the linear operator $L:\mathcal{D}(L) \to \mathcal{E}$  that satisfies
\begin{equation}
\label{eq:Lie_gen}
L \zeta(u) = \lim_{t \downarrow 0} \frac{K^t \zeta(u)-\zeta(u)}{t} \qquad \zeta \in \mathcal{D}(L)
\end{equation}
for all $u \in \overline{\mathcal{U}}$. \hfill $\diamond$
\end{definition}
\begin{remark}[Infinitesimal generator]
The Lie generator is reminiscent of the infinitesimal generator of strongly continuous semigroups (\cite{Engel_Nagel}). However, the limit in \eqref{eq:Lie_gen} is defined pointwise, while it is defined in the strong sense in the case of the infinitesimal generator. In fact, the infinitesimal generator of the Koopman semigroup cannot be defined unless the Koopman semigroup is strongly continuous (i.e. $\lim_{t \downarrow 0} \|K^t \zeta - \zeta\|=0$ $\forall \zeta \in \mathcal{E}$). This property does not hold in our setting, but can be satisfied with the mixed topology on $C(\overline{\mathcal{U}})$, as shown in \cite{Farkas}.
\end{remark}

The Lie generator enjoys a few properties (e.g. dense domain, bounded resolvent) and we refer to \cite{Neuberger} for more details. In the case of a Koopman semigroup associated with a semiflow generated by the abstract differential equation \eqref{eq:abstract_diff}, it follows from the chain rule property that the Lie generator is given by
\begin{equation*}
L \zeta(u) = \lim_{t \downarrow 0} \frac{\zeta(\varphi^t(u))-\zeta(u)}{t} = D_{D_t \varphi^t(u)} \zeta(u) = D_{W(u)} \zeta(u)
\end{equation*}
where $D_{W(u)} \zeta(u)$ denotes the Gâteaux derivative of $\zeta$ at $u$ in the direction $W(u)$:
\begin{equation*}
D_{W(u)} \zeta(u) = \lim_{\lambda \rightarrow 0} \frac{\zeta(u+\lambda W(u)) - \zeta(u)}{\lambda} \,.
\end{equation*}
Note that this can be interpreted as the Lie derivative associated with the infinite-dimensional vector field $W(\cdot)$.

\subsection{Finite-dimensional representation}

It is convenient to approximate the Koopman semigroup $K^t$ or the Lie generator $L$ in a finite-dimensional subspace of $\mathcal{E}$. Toward that end, we can consider the compressions $K_n^t=P_n K^t|_{\mathcal{E}_n}$ and $L_n=P_n L|_{\mathcal{E}_n}$, where $\mathcal{E}_n \subset \mathcal{D}(L)$ is a $n$-dimensional subspace of $\mathcal{E}$ and $P_n:\mathcal{E} \to \mathcal{E}_n$ is a projection operator. Suppose that $\mathcal{E}_n$ is spanned by the basis of functionals $\{\zeta_i\}_{i=1}^n$. In this basis, the finite-dimensional operators $K^t_n$ and $L_n$ can be represented by the matrices $\ve{K}\in \mathbb{R}^{n \times n}$ and $\ve{L}\in \mathbb{R}^{n \times n}$, respectively, which are defined so that
\begin{equation}
\label{eq:finite_approx}
K^t_n \zeta_j = \sum_{i=1}^n \ve{K}_{ij} \zeta_i \quad \textrm{and} \quad L_n \zeta_j = \sum_{i=1}^n \ve{L}_{ij} \zeta_i \,.
\end{equation}

The choice of the basis functions is crucial, as it affects the accuracy of the approximation and the performance of the methods based on this approximation (see below). However, finding the optimal set of basis functions is not trivial, and may require \textit{a priori} knowledge on the system.

\section{Spectral analysis and Extended Dynamic Mode Decomposition}
\label{sec:spectral}

The spectral properties of the Koopman operator reveal important geometric properties of the underlying dynamics (see e.g. \cite{Mezic}, and \cite{Mauroy_Mezic,Nakao_SICE} in the context of phase reduction). In this section, we exploit the proposed framework  for infinite-dimensional systems and compute the spectrum of the Koopman operator from data. This yields a generalization of the Extended Dynamical Mode Decomposition method for infinite-dimensional systems.

\subsection{Spectrum of the Koopman operator}
\label{sec:Koopman_spectrum}

We consider the spectrum of the Lie generator \eqref{eq:Lie_gen}, i.e. the set of (Koopman) eigenvalues $\lambda$ such that
$L \zeta_{\lambda} = \lambda \, \zeta_{\lambda}$ for some (Koopman) eigenfunctional $\zeta_{\lambda} \in \mathcal{E}$.

\paragraph{Case of linear systems.} It is well-known that the spectrum of linear finite-dimensional systems is contained in the spectrum of the related Koopman operator. As shown in \cite{Mezic_book} and \cite{Nakao_PDE}, this result also holds for infinite-dimensional systems. Consider a linear system $\dot{u} = A u$, $u\in \mathcal{U}$, and suppose that $\lambda$ is an eigenvalue of $A$, so that there exists an eigenfunction $w_\lambda\in \mathcal{D}(A^*)$ with $A^* w_\lambda =\bar{\lambda} w_\lambda$, where $A^*$ denotes the adjoint operator of $A$ and $\bar{\lambda}$ is the complex conjugate of $\lambda$. Then the functional $\zeta_{\lambda}(\cdot) = \langle \cdot, w_\lambda \rangle$ satisfies
\begin{equation*}
\begin{split}
L \xi_{\lambda}(u) = D_{Au} \langle u, w_\lambda \rangle= \langle Au, w_\lambda \rangle = \langle u, A^* w_\lambda \rangle & = \langle u, \bar{\lambda} w_\lambda \rangle \\
&  = \lambda \, \xi_{\lambda}(u) \,,
\end{split}
\end{equation*}
so that $\lambda$ is a Koopman eigenvalue. Moreover, it is easy to verify that $(\zeta_{\lambda})^\alpha$ is a Koopman eigenfunctional associated with the Koopman eigenvalue $\alpha \lambda$, provided that it belongs to the space $\mathcal{E}$. If $A$ generates a strongly continuous (linear) semigroup $(\varphi^t)_{t \geq 0}=(T^t)_{t \geq 0}$, then the spectral mapping theorem implies that $e^{\lambda t}$ is an eigenvalue of the operator $T^t$ (see e.g. \cite{Engel_Nagel}, Chapter IV). It follows that 
\begin{equation*}
K^t \zeta_{\lambda} = \langle T^t u, w_\lambda \rangle = \langle u, (T^t)^* w_\lambda \rangle = \langle u, e^{\bar{\lambda}t} w_\lambda \rangle = e^{\lambda t} \, \zeta_{\lambda}
\end{equation*}
and $e^{\lambda t}$ is also an eigenvalue of $K^t$.

\subsection{Extended Dynamic Mode Decomposition for infinite- dimensional systems}
\label{sec:EDMD}

Extended Dynamic Mode Decomposition (EDMD) is a data-driven method that builds a finite-dimensional approximation of the Koopman semigroup and computes the approximate spectral properties of the operator (\cite{Rowley_EDMD}). It can be easily extended to the infinite-dimensional framework that we consider here. 

Suppose we have access to a set of $m$ pairs $(u_k,\varphi^{t_s}(u_k)) \in \overline{\mathcal{U}}^2$, where $t_s$ is a sampling time, in such a way that we can measure the values of $n$ functionals $(\zeta_i(u_k), \zeta_i(\varphi^{t_s}(u_k))) \in \mathbb{R}^2$ for all $i\in \{1,\dots,n\}$. The generalized EDMD method proceeds with the following steps.
\begin{enumerate}
\item[1.] Compute the data matrices
\begin{equation}
\label{eq:data1}
\ve{\Theta}_1 = \left(\begin{array}{ccc}
\zeta_1(u_1) & \cdots & \zeta_n(u_1) \\
\vdots &  & \vdots \\
\zeta_1(u_m) & \cdots & \zeta_n(u_m) \end{array} \right)
\end{equation}
and
\begin{equation}
\label{eq:data2}
\ve{\Theta}_2 = \left(\begin{array}{ccc}
\zeta_1(\varphi^{t_s}(u_1)) & \cdots & \zeta_n(\varphi^{t_s}(u_1)) \\
\vdots &  & \vdots \\
\zeta_1(\varphi^{t_s}(u_m)) & \cdots & \zeta_n(\varphi^{t_s}(u_m)) \end{array} \right) \,.
\end{equation}
\item[2.] Provided that $m\geq n$, a matrix approximation of \mbox{$K^{t_s}_n=P_n K^{t_s}|_{\mathcal{E}_n}$} is given by the least squares solution $\ve{K} = \ve{\Theta}_1^\dagger \ve{\Theta}_2$, where $\ve{\Theta}^\dagger$ denotes the Moore-Penrose pseudoinverse of $\ve{\Theta}$. Note that, in this case, $P_n$ is the discrete orthogonal projection
\begin{equation}
\label{eq:disc_proj}
P_n \zeta = \underset{\tilde{\zeta} \in \mathcal{E}_n}{\textrm{argmin}} \sum_{k=1}^m |\zeta(u_k)-\tilde{\zeta}(u_k)|^2 \,.
\end{equation}
\item[3.] The eigenvalues $\lambda_K$ of $K^{t_s}$ are approximated by the eigenvalues of $\ve{K}$ and estimates of the eigenvalues $\lambda_L$ of the Lie generator are given by $\lambda_L=\log(\lambda_K)/t_s$. Moreover, Koopman eigenfunctionals are approximated in the basis of functionals by the components of the corresponding (right) eigenvectors of $\ve{K}$.
\end{enumerate}

The method only requires to know the samples $u_k$ in a weak sense, i.e. through the values $\zeta_j(u_k)$ of a finite number of functionals. For instance, these values could be the weighted averages of $u_k:X \to \mathbb{R}$ over the definition domain $X$. For the specific choice of evaluation functionals $\zeta_j(u)=u(x_j)$ with $x_j\in X$, we recover the classical DMD method (\cite{Tu_DMD}) applied to a discretized version of the infinite-dimensional system \eqref{eq:abstract_diff}.

\begin{remark}
The method is more general than the standard EDMD method in that it only requires to know the samples $u_k$ in a weak sense, i.e. through the values $\zeta_j(u_k)$ of a finite number of functionals. For instance, these values could be the weighted averages of $u_k:X \to \mathbb{R}$ over the definition domain $X$. For the specific choice of evaluation functionals $\zeta_j(u)=u(x_j)$ with the sample points $x_j\in X$, we recover the classical DMD method (\cite{Tu_DMD}) applied to a spatially discretized version of the infinite-dimensional system \eqref{eq:abstract_diff}. Similarly, evaluation functionals of the form $\zeta_j(u)=\psi(u(x_j))$, where $\psi$ is a basis function, yield the EDMD method \cite{Rowley_EDMD} applied to the discretized infinite-dimensional system.
\end{remark}

\begin{remark}[Convergence properties]\label{rem:converg} The EDMD method must be used with some care since it might yield spurious eigenvalues. Similarly to \cite{Korda_convergence}, convergence properties should be characterized as $m,n\rightarrow \infty$ and the general validity of the spectral mapping theorem should be investigated. We leave these questions for future research.
\end{remark}

\subsection{Numerical example}

We illustrate the generalized EDMD method with the Burgers equation
\begin{equation}
\label{eq:dyn_Burgers}
\dot{u} = - u \frac{\partial u}{\partial x} + \frac{\partial^2 u}{\partial x^2} \qquad u \in L^2([-1,1])
\end{equation}
associated with homogeneous Dirichlet boundary conditions $u(-1)=u(1)=0$. The dynamics \eqref{eq:dyn_Burgers} is conjugated to the linear diffusion dynamics $\dot{v} = \partial^2 v/\partial t^2$ through the so-called Cole-Hopf transformation (\cite{Hopf}). As explained in \cite{Nakao_PDE,Page}, this implies that the Koopman spectrum associated with the Burgers dynamics coincides with the Koopman spectrum associated with the linear diffusion dynamics. It follows that Koopman eigenfunctionals are related to eigenfunctions of the diffusion operator (see Section \ref{sec:Koopman_spectrum}) and, in particular, they are of the form $\zeta(v)=\langle sin(k\pi x/2), v(x) \rangle^\alpha$ (in the new variable $v$). The associated Koopman eigenvalues are given by $\lambda=\alpha (k\pi/2)^2$.

\begin{figure}[h]
\begin{center}
\includegraphics[width=8.4cm]{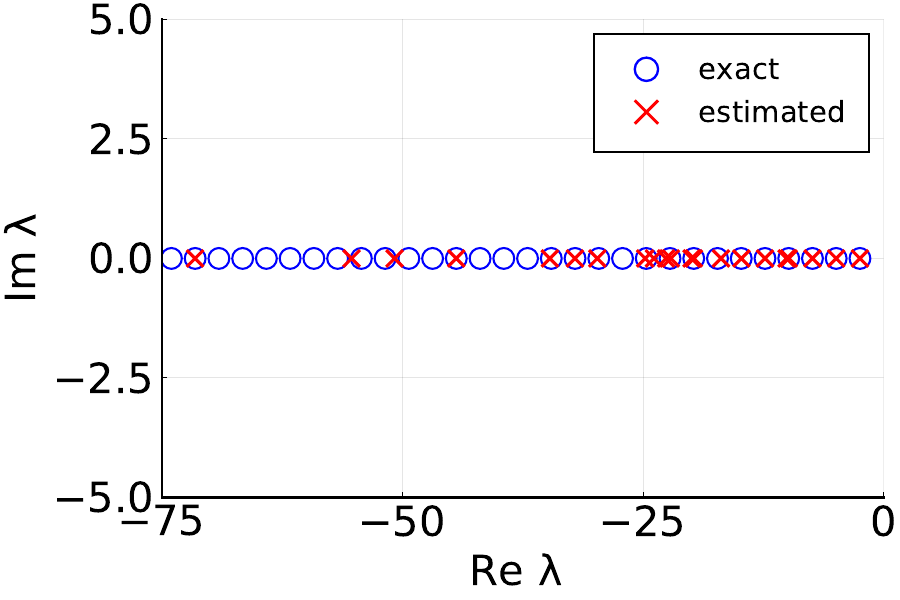}    
\caption{The generalized EDMD method is used to compute the Koopman eigenvalues associated with the Burgers dynamics \eqref{eq:dyn_Burgers}. Dominant eigenvalues $-\alpha (\pi/2)^2$, $\alpha\in \mathbb{N}$ (blue circles) are correctly estimated (red crosses).} 
\label{fig:spectrum_Burgers}
\end{center}
\end{figure}

We use the dynamics \eqref{eq:dyn_Burgers} to generate $m=50$ data-pairs taken from $10$ trajectories (with random, arbitrarily chosen initial conditions of the form $(x^2-1)\cos(a\pi x+b \pi)$, $a,b\in[0,1]$). The sampling time is $t_s=0.2$. Estimates of the Koopman eigenvalues are computed through our generalized EDMD method, with $n=27$ basis functionals $\zeta_{j,k,l}=\langle \cos(a_j(\pi x/2)+b_j\pi/2), (u(x))^k \rangle^l$, with $(j,k,l) \in \{1,2,3\}^3$ and where $a_j,b_j$ are randomly chosen over the interval $[0,1]$. As shown in Fig. \ref{fig:spectrum_Burgers}, dominant eigenvalues $-\alpha (\pi/2)^2$, $\alpha \in \mathbb{N}$, are captured. This is consistent with the results presented in \cite{Page} and shows the importance of selecting an augmented basis of nonlinear functionals to capture other eigenvalues than the principal ones of the form $(k \pi/2)^2$, $k\in \mathbb{N}$. Note that other (complex) eigenvalues may appear over different numerical tests. They should be considered in light of future theoretical analysis (see Remark \ref{rem:converg}), and more advanced methods should be proposed to distinguish true eigenvalues from spurious ones.

\section{Identification of infinite-dimensional systems}
\label{sec:identif}

In this section, we use the Koopman operator framework for infinite-dimensional systems in the context of identification. Our goal is to identify the coefficients $c_i \in \mathbb{R}$ of the infinite-dimensional dynamics
\begin{equation}
\label{eq:abstract_diff2}
\dot{u} = W(u) = \sum_i c_i W_i(u) \,,
\end{equation}
with $u \in \mathcal{U}=L^2(X)$ (where $X\subset \mathbb{R}^p$ is a compact set), using $m$ data pairs $(u_k,\varphi^{t_s}(u_k))$ generated by the dynamics \eqref{eq:abstract_diff2}. We assume that the operators $W_i:\mathcal{D}(W_i) \to \mathcal{U}$ are known \textit{a priori}, so that this can be seen as a parameter estimation problem. Note that \eqref{eq:abstract_diff2} may be described by a partial differential equation, 
but the proposed method is not limited to that case (see Section \ref{sec:identif_num}). 

\subsection{Lifting identification method}

The lifting identification method proposed in our previous work (\cite{Mauroy_Goncalves2}) is generalized to the case of infinite-dimensional systems. This method consists in three steps.
\begin{enumerate}
\item[1.] \emph{Lifting of the data.} We compute the data matrices \eqref{eq:data1} and \eqref{eq:data2} with the basis of functionals
\begin{equation}
\label{eq:basis_ident}
\zeta_i(u) = \langle W_i(u),w \rangle
\end{equation}
where $\langle \cdot, \cdot \rangle$ denotes the inner product in $L^2(X)$ and $w \in L^2(X)$ is a weighting function. We suppose without loss of generality that $W_1(u)=u$ (possibly with the coefficient $c_1=0$ in \eqref{eq:abstract_diff2}), so that the linear functional $\zeta_1(u)=\langle u,w \rangle$ belongs to $\mathcal{E}_n$.
\smallskip
\item[2.] \emph{Identification of the Lie generator.} We compute the matrix representation $\ve{K}=\ve{\Theta}_1^\dagger \ve{\Theta}_2$ of the compression $K_n^{t_s}$ in the subspace $\mathcal{E}_n=\textrm{span}(\zeta_1,\dots,\zeta_n)$ (step 2 in Section \ref{sec:EDMD}). Then we obtain a finite-dimensional approximation $\widetilde{\ve{L}}^{(t_s)}$ of the Lie generator by taking the matrix logarithm
\begin{equation}
\label{eq:L_tilde}
\widetilde{\ve{L}}^{(t_s)} = \frac{1}{t_s} \log (\ve{\Theta}_1^\dagger \ve{\Theta}_2) \,.
\end{equation}
Note that this approximation is not equal to the matrix representation $\ve{L}$ of the $L_n$ (see \eqref{eq:finite_approx}).
\smallskip
\item[3.] \emph{Identification of the coefficients.} Estimates $\hat{c}_i$ of the coefficients $c_i$ are given by the entries of the first column of $\widetilde{\ve{L}}^{(t_s)}$, i.e. $\hat{c}_i = \widetilde{\ve{L}}^{(t_s)}_{i1}$.
\end{enumerate}
\begin{remark} For the specific choice of basis functionals of the form $\zeta_j(u)=\psi(u(x_j))$, with $x_j\in X$ and where $\psi$ is a basis function, one recovers the original lifting identification method \cite{Mauroy_Goncalves2} applied to a spatially discretized version of the infinite-dimensional system. The proposed method is more general since it allows any basis functionals.
\end{remark}

The lifting identification method does not require to compute time derivatives and is therefore an alternative to direct methods for PDE identification such as those proposed in \cite{Kutz_PDE, PDE_ident_neural_net, Yuan_PDE, Grigoriev}. It is also noticeable that the use of basis functionals of the form \eqref{eq:basis_ident} bears similarity to the method developed in \cite{Grigoriev}, which makes a clever use of a weak formulation of the data in space and time. We note that this method requires sufficiently long time-series to allow accurate time integration, while our method can deal with data pairs belonging to different trajectories.

\subsection{Convergence results}

Now we show that the estimated coefficients converge to the true coefficients as $t_s \rightarrow 0$. This is summarized in the following proposition.
\begin{proposition}
\label{prop1}
Let $\varphi^t$ be the continuous flow generated by \eqref{eq:abstract_diff2} and let $\mathcal{E}_n \subset \mathcal{D}(L) $ be the subspace spanned by the basis functionals $\zeta_i(\cdot)=\langle W_i(\cdot),w \rangle \in C(\overline{\mathcal{U}})$, with $\zeta_1(\cdot)=\langle \cdot,w \rangle$. Assume that the pairs $(u_k,\varphi^{t_s}(u_k)) \in \overline{\mathcal{U}}^2$ are such that the $n \leq m$ vectors $(\zeta_i(u_1) , \dots , \zeta_i(u_m))$, $i=1,\dots,n$, are linearly independent. Then,
\begin{equation*}
\lim_{t_s \rightarrow 0} \widetilde{\ve{L}}^{(t_s)}_{i1} = c_i
\end{equation*}
where $\widetilde{\ve{L}}^{(t_s)}$ is given by \eqref{eq:L_tilde}.
\end{proposition}
\begin{pf}
The discrete orthogonal projection $P_n$ \eqref{eq:disc_proj} is well-defined since the vectors \\
$(\zeta_i(u_1) , \dots , \zeta_i(u_m))$, $i=1,\dots,n$, are linearly independent. It is clear that
\begin{equation*}
\frac{1}{t_s} \log K_n^{t_s} \zeta_1(u) = \sum_{i=1}^n \widetilde{\ve{L}}_{i1}^{(t_s)} \zeta_i(u)
\end{equation*}
and we also have
\begin{equation}
\label{eq:L_xi_1}
L \zeta_1(u)  = D_{W(u)} \langle u,w \rangle = \langle W(u),w \rangle = \sum_{i=1}^n c_i \zeta_i(u) \,.
\end{equation}
Since $P_n \zeta_i = \zeta_i$ for all $i=1,\dots,n$, it follows that
\begin{equation*}
\left|\sum_{i=1}^n \left(\widetilde{\ve{L}}^{(t_s)}_{i1} - c_i \right) \zeta_i(u)\right| =  \left| \left(\frac{1}{t_s} \log K^{t_s}_n  - P_n L P_n \right) \zeta_1(u) \right|.
\end{equation*}
For $t_s$ small enough, one has $\log A_n^{t_s}=P_n L P_n$, where $A^{t_s}_n=\mathcal{E}_n \to \mathcal{E}_n$ is the finite-dimensional operator $A^{t_s}_n:=e^{t_s P_n L P_n}$, so that
\begin{equation}
\label{eq:L_c_i}
\left|\sum_{i=1}^n \left(\widetilde{\ve{L}}^{(t_s)}_{i1} - c_i \right) \zeta_i(u)\right| = \frac{1}{t_s} \left| \left(  \log K^{t_s}_n - \log A_n^{t_s} \right) \zeta_1(u) \right|.
\end{equation}
Since the basis functionals $\zeta_i$ are linearly independent,
\begin{equation}
\label{eq:lim_zero}
\lim_{t_s \rightarrow 0}\frac{1}{t_s}|(\log A_n^{t_s} - \log K_n^{t_s} ) \zeta_1(u) | = 0 \qquad \forall u
\end{equation}
implies that
\begin{equation*}
\lim_{t_s \rightarrow 0} \left|\widetilde{\ve{L}}^{(t_s)}_{i1} - c_i\right|=0 \,.
\end{equation*}
Thus it remains to show that \eqref{eq:lim_zero} holds.

Since \mbox{$\lim_{t_s \rightarrow 0}\|A^{t_s}_n-I\|=0$}, it is easy to check that
\begin{equation}
\label{eq:term1}
\lim_{t_s \rightarrow 0} \frac{1}{t_s}\left\|\log A_n^{t_s} - (A_n^{t_s}-I) \right\| = 0
\end{equation}
[note that $\lim_{t \rightarrow 0} (\log f(t)-(f(t)-1))/t=0$ for all $f \in C^1$ such that $\lim_{t \rightarrow 0} f(t)=1$] and it follows similarly that
\begin{equation}
\label{eq:term2}
\lim_{t_s \rightarrow 0} \frac{1}{t_s}\left\|\log K_n^{t_s} - (K_n^{t_s}-I) \right\| = 0\,.
\end{equation}
Since the mapping $t \mapsto (A_n^{t} - K_n^{t}) \zeta(u)$ is differentiable, the mean value theorem implies that
\begin{equation*}
\frac{1}{t_s}\left|(A_n^{t_s} - K_n^{t_s} ) \zeta(u) \right| = \left|\frac{d}{dt} \left(A^t_n K^{t_s-t}\zeta(u)\right)|_{t=\tau} \right|
\end{equation*}
for all $\zeta \in \mathcal{E}_n$ and $u \in \overline{\mathcal{U}}$, and for some $\tau \in [0,t_s]$. Then we can write
\begin{equation*}
\begin{split}
& \frac{1}{t_s}\left|(A_n^{t_s} - K_n^{t_s} ) \zeta(u) \right| \\
& = \left | A_n^\tau (P_n L P_n-P_n L) K^{t_s-\tau} \zeta(u) \right| \\
& \leq | A_n^\tau P_n L P_n (K^{t_s-\tau} \zeta(u) - \zeta(u)) | \\
& \quad + |A_n^\tau P_n L ( P_n \zeta(u) - \zeta(u)) | + | A_n^\tau P_n L (\zeta(u) - K^{t_s-\tau} \zeta(u)) | \\
& \leq | A_n^\tau P_n L P_n (\zeta(\varphi^{t_s-\tau}(u)) - \zeta(u)) | \\
& \quad + | A_n^\tau P_n L (\zeta(u) - \zeta(\varphi^{t_s-\tau}(u))) |
\end{split}
\end{equation*}
where we used the fact that $P_n\zeta = \zeta$ and that $K^t$ and $L$ commute. Since the functionals $A_n^\tau P_n L P_n \zeta$ and $A_n^\tau P_n L \zeta$ are continuous and the flow $\varphi^{t}$ is continuous in $t$, we have
\begin{equation*}
\begin{split}
\lim_{t_s \rightarrow 0} | A_n^\tau P_n L P_n (\zeta(\varphi^{t_s-\tau}(u)) - \zeta(u)) | =0 \\
\lim_{t_s \rightarrow 0} | A_n^\tau P_n L ( \zeta(u) - \zeta(\varphi^{t_s-\tau}(u))) | = 0
\end{split}
\end{equation*}
so that
\begin{equation}
\label{eq:term3}
\lim_{t_s \rightarrow 0} \frac{1}{t_s}\left|(A_n^{t_s} - K_n^{t_s} ) \zeta(u) \right| = 0 \,.
\end{equation}
Finally \eqref{eq:term1}, \eqref{eq:term2} and \eqref{eq:term3} yield
\begin{equation*}
\begin{split}
& \lim_{t_s \rightarrow 0}\frac{1}{t_s}|(\log A_n^{t_s} - \log K_n^{t_s} ) \zeta(u) | \\
& \, \leq \lim_{t_s \rightarrow 0} \frac{1}{t_s}\left( |(\log A_n^{t_s} - (A_n^{t_s}-I) ) \zeta(u) | \right. \\
& \quad \left. + |(\log K_n^{t_s} - (K_n^{t_s}-I) ) \zeta(u) | + |(A_n^{t_s} - K_n^{t_s} ) \zeta(u) | \right) \\
& \,  = 0 \,.
\end{split}
\end{equation*}
This concludes the proof.
\end{pf}
Although the result requires that the sampling time $t_s$ tend to zero, accurate results can be obtained even if $t_s$ is not so small (see Section \ref{sec:identif_num}). However, large sampling times might affect the method performance. 
In particular, if the sampling time is too large, the principal branch of the matrix logarithm $1/t_s \log A^{t_s}_n = 1/t_s \log e^{t_s P_n L P_n}$ is not equal to $P_n L P_n$. This is related to the so-called system aliasing issue \cite{Yue} (see also \cite{Mauroy_Goncalves2} for more details).

When an arbitrary basis of functionals is chosen so that $\langle W(\cdot),w \rangle \notin \textrm{span} \{\zeta_1,\dots,\zeta_n\}$, we have the following corollary.
\begin{corollary}
\label{col1}
Under the assumptions of Proposition \ref{prop1}, but with $\zeta_i(\cdot)=\langle V_i(\cdot),w \rangle$, $V_i:\mathcal{D}(V_i) \to \mathcal{U}$ (and $V_1=I$), the estimated nonlinear operator $\hat{W}^{(t_s)} = \sum_{i=1}^n \widetilde{\ve{L}}^{(t_s)}_{i1} V_i$ satisfies
\begin{equation*}
\lim_{t_s \rightarrow 0} \hat{W}^{(t_s)} = \underset{\tilde{W} \in \textrm{span}(V_1,\dots,V_n)}{\textrm{argmin}} \sum_{k=1}^m \left|\left\langle \tilde{W}(u_k)-W(u_k) ,w\right\rangle \right|^2 \,.
\end{equation*}
\end{corollary}
\begin{pf}
It follows from \eqref{eq:L_xi_1}, \eqref{eq:L_c_i} and \eqref{eq:lim_zero} that
\begin{equation*}
\lim_{t_s \rightarrow 0} \sum_{i=1}^n \widetilde{\ve{L}}^{(t_s)}_{i1} \zeta_i(\cdot) = P_n \langle W(\cdot),w \rangle \,,
\end{equation*}
which implies that
\begin{equation*}
\lim_{t_s \rightarrow 0} \left \langle \hat{W}^{(t_s)}(\cdot),w \right\rangle =  P_n \langle W(\cdot),w \rangle \,.
\end{equation*}
The result follows from the definition of the discrete orthogonal projection \eqref{eq:disc_proj}.
\end{pf}
\begin{remark}
The proof of Proposition \ref{prop1} is adapted from a proof in \cite{Mauroy_Goncalves2}, but does not rely on semigroup theory. In particular, only weak (i.e. pointwise) convergence properties are used in the proof of Proposition \ref{prop1}, while strong convergence is considered in the proof in \cite{Mauroy_Goncalves2}. For this reason, strong convergence results can be obtained in \cite{Mauroy_Goncalves2} in the case of arbitrary bases, but a weaker result is obtained here (see Corollary \ref{col1}). In this context, considering a space of bounded continuous functionals equipped with a mixed topology would allow to exploit the strong continuity property of the semigroup, as shown in \cite{Farkas}, and possibly recover stronger convergence results.
\end{remark}

\subsection{Case of linearly dependent basis functionals}
\label{sec:dependent_fct}

The basis functionals \eqref{eq:basis_ident} might not be linearly independent, even if the operators $W_i$ are linearly independent (see Section \ref{sec:graphon}). In this case, the result of Proposition \ref{prop1} does not hold and in particular the estimated coefficients $\hat{c}_i\triangleq \widetilde{\ve{L}}^{(t_s)}_{i1}$ do not approximate the exact coefficients $c_i$. However, it follows from \eqref{eq:L_c_i}-\eqref{eq:lim_zero} that these coefficients satisfy the equality
\begin{equation*}
\sum_{i=1}^N (\hat{c}_i-c_i)  \zeta_i(u) = 0
\end{equation*}
as $t_s$ goes to zero. Considering several weighting functions $w^{(j)}$, with $j=1,\dots,J$, we can use the proposed identification method with $J$ different sets of basis functionals $\zeta_i^{(j)}=\langle W_i(\cdot),w^{(j)}\rangle$, yielding several sets of values $\hat{c}_i^{(j)}$ that satisfy the equations
\begin{equation*}
\sum_{i=1}^N (\hat{c}_i^{(j)}-c_i) \, \zeta^{(j)}_i(u) =  0 \qquad j=1,\dots,J.
\end{equation*}
Considering the above set of equations at $u=u_k$, for $k=1,\dots,m$, we obtain the matrix equality
\begin{equation}
\label{eq:equation_depend}
\ve{\Theta}_{full} \left(\begin{array}{c}
c_1  \\
 \vdots \\
 c_N
\end{array} \right) = \ve{b}
\end{equation} 
with
\begin{equation*}
\ve{\Theta}_{full} = \left(\begin{array}{c}
\ve{\Theta}_1^{(1)}  \\
 \vdots \\
 \ve{\Theta}_1^{(J)} 
\end{array} \right), \qquad \ve{b}= \left(\begin{array}{c}
\ve{\Theta}_1^{(1)} \left(\begin{array}{c}
\hat{c}_1^{(1)}  \\
 \vdots \\
 \hat{c}_N^{(1)}
\end{array} \right) \\
 \vdots \\
 \ve{\Theta}_1^{(J)} \left(\begin{array}{c}
 \hat{c}_1^{(J)}  \\
  \vdots \\
  \hat{c}_N^{(J)}
 \end{array} \right)
\end{array} \right),  
\end{equation*} 
and where $\ve{\Theta}_1^{(j)}$ is the data matrix \eqref{eq:data1} obtained with basis functionals $\zeta_i^{(j)}$. Provided that $J$ is large enough and the choice of weighting functions $w^{(j)}$ is appropriate, the matrix $\ve{\Theta}_{full}$ can be full rank so that \eqref{eq:equation_depend} admits a unique solution $\ve{\Theta}_{full}^\dagger \ve{b}$. In this case, the coefficients $c_i$ are recovered despite the fact that every set of basis functionals is linearly dependent.

\begin{remark}
\label{rem:several_weights}
Several weighting functions can also be used when the basis functionals are linearly independent (provided that the observations are available in practice). Then, the sets of estimated coefficients $\hat{c}_i^{(j)}$ can be averaged to improve the accuracy of the results. Inconsistent results can also be discarded by comparing the different sets of coefficients.
\end{remark}

\subsection{Numerical examples}
\label{sec:identif_num}

We can now use the lifting identification method with two illustrating examples: a nonlinear partial differential equation and a nonlinear diffusive dynamics on a graphon.

\subsubsection{Nonlinear partial differential equation}
\label{sec:PDE_identif}

\begin{figure}[h]
\begin{center}
\includegraphics[width=7cm]{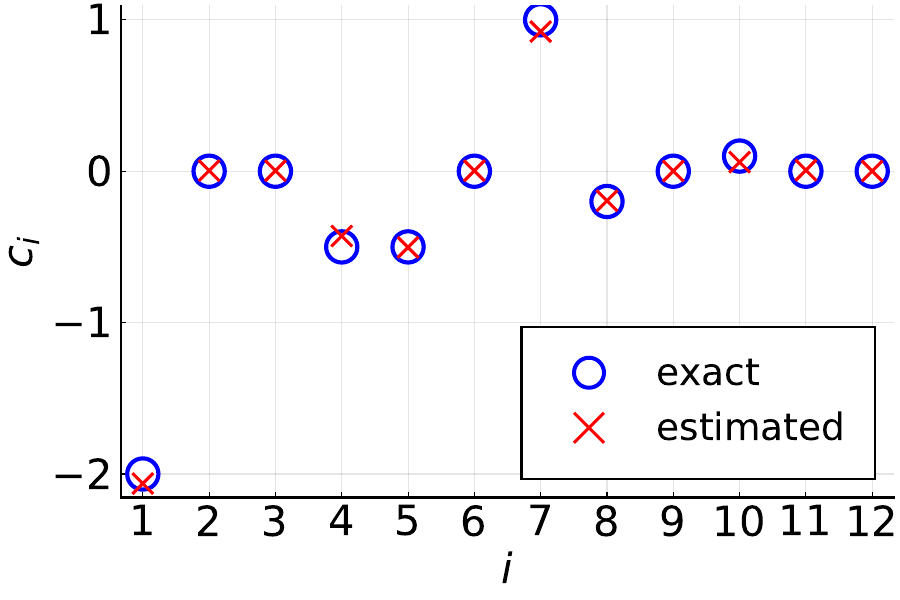}    
\caption{The lifting identification method is used to recover the PDE \eqref{eq:PDE_dyn}. The coefficients $c_i$ (blue circles) are estimated with small error (red crosses). Nonzero coefficients correspond to $W_1(u)=u$, $W_4(u)=\partial u/\partial x$,  $W_5(u)=u \, \partial u/\partial x$, $W_7(u) = \partial^2 u/\partial x^2$, $W_8(u) = u\, \partial^2 u/\partial x^2$, and $W_{10}(u)=\partial^3 u/\partial x^3$.} 
\label{fig:ident_PDE}
\end{center}
\end{figure}

We aim at identifying the coefficients of the dynamics
\begin{equation}
\label{eq:PDE_dyn}
\dot{u} = -2u - 0.5 (1+u) \frac{\partial u}{\partial x} + (1-0.2 u) \frac{\partial^2 u}{\partial x^2} + 0.1 \frac{\partial^3 u}{\partial x^3},
\end{equation}
$u\in L^2[0,5]$, with homogeneous Dirichlet boundary conditions $u(0)=u(5)=0$. The PDE is used to generate $m=50$ data pairs, taken from $25$ trajectories (with random, arbitrarily chosen initial conditions of the form $x(x-5)\cos(a\pi x/5+b \pi)$, $a,b\in[0,1]$). The sampling time is $t_s=0.3$. The lifting identification method is used with $n=12$ basis functionals \eqref{eq:basis_ident}, with the nonlinear operators
\begin{equation*}
W_i(u) \in \left\{ u^j \frac{\partial^k u}{\partial x^k},\, j\in\{0,1,2\},\, k\in\{0,1,2,3\} \right\}
\end{equation*}
and the weighting function $w(x)=\exp(-1/(1-(x/L)^2))$. Fig. \ref{fig:ident_PDE} shows that the coefficients $c_i$ are estimated with small error.

\subsubsection{Nonlinear diffusive dynamics on a graphon}
\label{sec:graphon} 

\begin{figure}[h]
\begin{center}
\includegraphics[width=7cm]{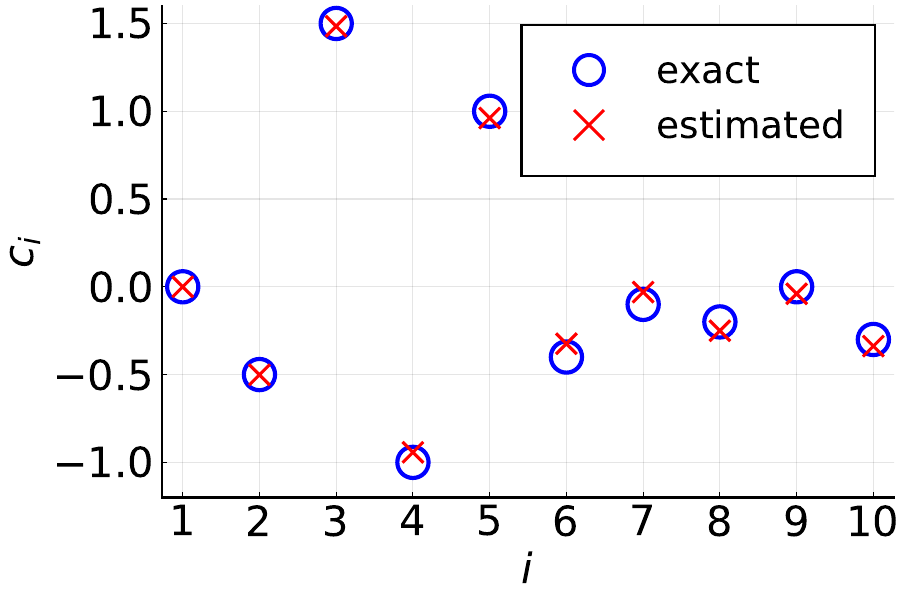}    
\caption{The lifting identification method is used to recover the dynamics \eqref{eq:graphon_dyn} on a graphon. The coefficients $c_i$ (blue circles) are correctly estimated (red crosses).}
\label{fig:ident_graphon}
\end{center}
\end{figure}

The graphon $G:[0,1]^2 \to [0,1]$ is used to describe the limit of a sequence of dense graphs, and can be interpreted as the infinite-dimensional version of an adjacency matrix. Here, we consider a nonlinear diffusive dynamics on the graphon $G(x,y)=1-0.4x-0.1y-0.2xy-0.3y^2$, which is described by the integro-differential equation
\begin{equation}
\label{eq:graphon_dyn}
\begin{split}
\dot{u}(t,x) & = -0.5u(t,x)+1.5u(t,x)^2-u(t,x)^3 \\
 & \quad + \int_0^1 G(x,y) \, (u(t,y)-u(t,x)) \, dy.
\end{split}
\end{equation}
We generate $m=60$ data pairs, taken from $30$ trajectories (with random initial conditions of the form $0.1 a \cos(b\pi x+b \pi)$, $a,b\in[0,1]$). The sampling time is $t_s=0.5$. The lifting identification method is used with $n=10$ basis functionals \eqref{eq:basis_ident}, with the operators
\begin{equation*}
\begin{split}
& W_1(u) = 1, \quad W_2(u)=u, \quad W_3(u)=u^2, \quad W_4(u)=u^3, \\
& W_5(u)(x) = \int_0^1 u(y)-u(x) \, dy, \\
& W_6(u)(x) = \int_0^1 x \, (u(y)-u(x)) \, dy,\\
& W_7(u)(x) = \int_0^1 y \, (u(y)-u(x)) \, dy, \\
& W_8(u)(x) = \int_0^1 xy \, (u(y)-u(x)) \, dy, \\
& W_9(u)(x) = \int_0^1 x^2 \, (u(y)-u(x)) \, dy, \\
& W_{10}(u)(x) = \int_0^1 y^2 \, (u(y)-u(x)) \, dy.
\end{split}
\end{equation*}
In particular, we can compute the basis functionals $\zeta_i(\cdot)=\langle W_i(\cdot),w \rangle$, with $i=5,\dots,8$ and we obtain
\begin{equation*}
\begin{split}
\zeta_5(u) & = \int_0^1 \int_0^1 w(x) (u(y)-u(x)) \,dx dy \\
& = C_1 \, \langle 1,u(x) \rangle + C_2  \, \langle w(x), u(x) \rangle \\
\zeta_6(u) & = \int_0^1 \int_0^1 w(x) \, x \, (u(y)-u(x)) \, dx dy \\
& = C_3 \, \langle 1 ,u(x) \rangle + C_2 \, \langle x \, w(x),u(x) \rangle \\
\zeta_7(u) & = \int_0^1 \int_0^1 w(x) \, y \, (u(y)-u(x)) \, dx dy \\
& = C_1 \, \langle x ,u(x) \rangle + C_4 \, \langle w(x),u(x) \rangle \\
\zeta_8(u) & = \int_0^1 \int_0^1 w(x) \, xy \, (u(y)-u(x)) \, dx dy \\
& = C_3 \, \langle x ,u(x) \rangle + C_4 \, \langle x \, w(x),u(x) \rangle
\end{split}
\end{equation*}
with $C_1=\int_0^1 w(x) dx$, $C_2=-1$, $C_3=\int_0^1 x \, w(x) dx$, and $C_4=-1/2$. Then it is easy to see that the functionals $\zeta_i$, with $i=5,\dots,8$, are linearly dependent for any weight function $w$. We therefore follow the procedure described in Section \ref{sec:dependent_fct}, using the weighting functions $w^{(j)}(x)=x^j$, with $j=1,\dots,4$. As shown in Fig. \ref{fig:ident_graphon}, the coefficients are correctly estimated and, in particular, the graphon is identified.

\subsubsection{Numerical performance}

In this section, we compare the lifting identification method with a simple direct identification method. For this latter method, the time derivative $\dot{u}_k$ at $u_k$ is estimated through (forward) finite differences
\begin{equation*}
\dot{u}_k = \frac{\varphi^{t_s}(u_k)-u_k}{t_s}
\end{equation*}
or equivalently
\begin{equation*}
\langle \dot{u}_k,w \rangle = \frac{\langle \varphi^{t_s}(u_k),w \rangle- \langle u_k,w \rangle}{t_s} = \frac{\zeta_1(\varphi^{t_s}(u_k))- \zeta_1(u_k)}{t_s} \,.
\end{equation*}
The estimated coefficients $\hat{c}_i$ are obtained through least squares regression of $\langle \dot{u},w \rangle$ over basis functionals of the form \eqref{eq:basis_ident}, i.e.
\begin{equation*}
\begin{pmatrix}
\hat{c}_1 \\
\vdots \\
\hat{c}_n
\end{pmatrix} = \ve{\Theta}_1^\dagger \begin{pmatrix}
\langle \dot{u}_1,w \rangle \\
\vdots \\
\langle \dot{u}_m,w \rangle
\end{pmatrix}
\end{equation*}
where $\ve{\Theta}_1$ is the data matrix \eqref{eq:data1}.\\

\begin{figure}[h]
\begin{center}
\includegraphics[width=7cm]{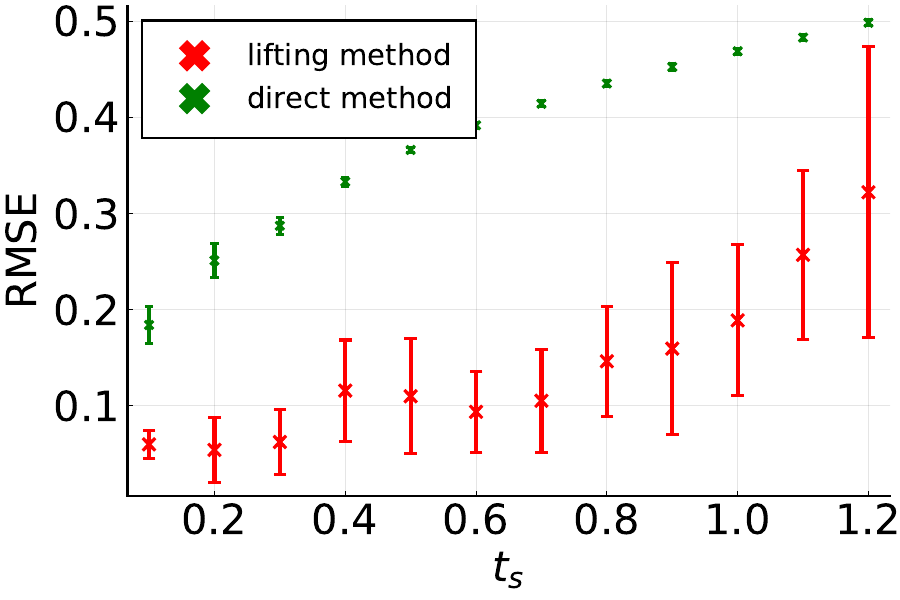}    
\caption{The lifting method is characterized by a smaller RMSE than a direct method based on least-squares regression of time derivatives. In particular, the error remains smaller for large sampling times. The dataset is generated with the dynamics \eqref{eq:PDE_dyn}. Crosses and error bars show respectively the mean and standard deviation of the RMSE (over $50$ experiments). Note that experiments characterized by $|\hat{c}_i^{(1)}-\hat{c}_i^{(2)}|>1$ for some $i$ are discarded.}
\label{fig:comparison}
\end{center}
\end{figure}

We apply the two methods on data generated by the PDE \eqref{eq:PDE_dyn}, with the same setting and set of parameters as in Section \ref{sec:PDE_identif}. However several values of the sampling time $t_s$ are considered and two weighting functions $w^{(1)}(x)=\exp(-1/(1-(x/L)^2))$ and $w^{(2)}(x)=\exp(-0.5/(1-(x/L)^2))$ are also used (see Remark \ref{rem:several_weights}). For each method, we compute two sets of coefficients $\hat{c}_i^{(1)}$ and $\hat{c}_i^{(2)}$ with the two weighting functions and take the mean $\hat{c}_i=(\hat{c}_i^{(1)}+\hat{c}_i^{(2)})/2$. We finally compute the root mean square error (RMSE) 
\begin{equation*}
RMSE=\frac{\sqrt{\sum_{i=1}^n (c_i-\hat{c}_i)^2}}{\sqrt{n}},
\end{equation*}
which is shown in Figure \ref{fig:comparison} for both methods as a function of the sampling time. We can see that the lifting method outperforms the direct method and, in particular, is characterized by a small RMSE even for large sampling times. However, it is characterized by a higher variability.

\begin{figure}[h]
\begin{center}
\includegraphics[width=7cm]{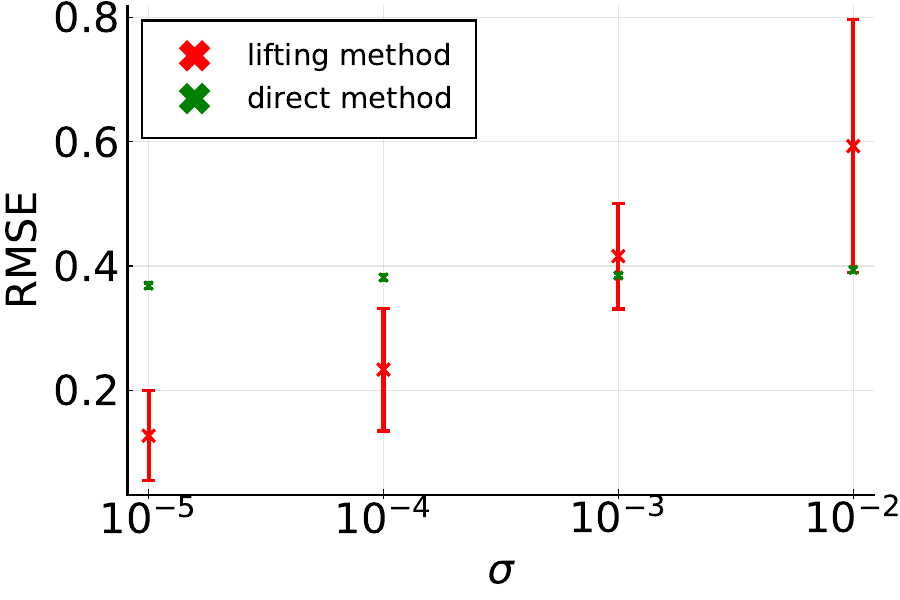}    
\caption{The lifting method is characterized by a smaller RMSE for low noise levels. However, the direct method is more robust to noise and outperforms the lifting method for larger noise levels. The dataset is generated with the dynamics \eqref{eq:PDE_dyn}. Crosses and error bars show respectively the mean and standard deviation of the RMSE (over $20$ experiments). Note that experiments characterized by $|\hat{c}_i^{(1)}-\hat{c}_i^{(2)}|>1$ for some $i$ are discarded.}
\label{fig:noise}
\end{center}
\end{figure}

We have also considered the effect of measurement noise. To do so, we have considered the same setting as above (with $t_s=0.5$) and we have added to the data a Gaussian noise with zero mean and standard deviation $\sigma \cdot \textrm{std}(\textrm{data})$, where $\textrm{std}(\textrm{data})$ stands for the standard deviation of the data. As shown in Figure \ref{fig:noise}, the lifting method still outperforms the direct method for low noise level (i.e. $0.001\%$ and $0.01\%$), but is not robust to larger noise level in which case the direct method yields better results. The variability of the results is also higher with the lifting method. This can be explained by the fact that our proposed parameter estimation method is biased and not consistent, due to the lifting of (noisy) data. In future work, noise robustness of the method should be improved. For instance, the weak formulation of the method could be further exploited to enhance noise robustness (see e.g. \cite{Grigoriev}).

\section{Conclusion}
\label{sec:conclu}

In this paper, the Koopman operator framework has been leveraged in the case of infinite-dimensional nonlinear systems. Building on previous theoretical works, we have considered a semigroup of composition operators and the associated Lie generator on a space of continuous functionals. A finite-dimensional representation of these operators has been proposed and used in the context of spectral analysis. This approach yields a generalization of the data-driven EDMD method for infinite-dimensional systems. We have also developed a novel identification method, which allows to identify nonlinear PDEs although it relies solely on linear techniques. This method has been complemented with convergence results.

The Koopman operator framework for infinite-dimensional systems is still in its infancy. Convergence properties of the finite-dimensional approximation of the Koopman operator should be thoroughly studied, in particular in light of the recent results by \cite{Farkas} in semigroup theory. This could provide some insight into the results obtained with the generalized EDMD method. The possibility to consider several weights functions should be further investigated and exploited. In particular, some guidelines to carefully select the weight functions could be provided. Finally, robustness to noise should be improved.





\end{document}